\documentclass[reqno]{amsart}
\usepackage{amsmath, amsthm, amsfonts, amsopn, graphicx}

    \theoremstyle{plain}
    \newtheorem{Thm}{Theorem}[section]
    \newtheorem{Prop}[Thm]{Proposition}
    
    \newtheorem*{Lemma*}{Lemma}

    \theoremstyle{definition}
    
    \newtheorem*{Def*}{Definition}
    \newtheorem{Example}[Thm]{Example}

    \theoremstyle{remark}
    \newtheorem{Remark}{Remark}
    \newtheorem*{Remark*}{Remark}

    \numberwithin{equation}{section}

    \newcommand{\ie}{\textit{i.e.}}

    \newcommand{\field}[1]{\mathbb{#1}}
    
    \newcommand{\Z}{\field{Z}}

    \newcommand{\C}{\field{C}}
    \DeclareMathSymbol{\fieldk}{\mathalpha}{AMSb}{"7C} 
    \newcommand{\abs}[1]{\lvert#1\rvert}
    \newcommand{\Abs}[1]{\left\lvert#1\right\rvert}
    
    \newcommand{\inner}[2]{\left\langle#1,#2\right\rangle}
    \DeclareMathOperator{\Tr}{Tr}
    \DeclareMathOperator{\Det}{Det}
    
    \DeclareMathSymbol{\normal}{\mathord}{AMSa}{"43}
    \DeclareMathOperator{\Log}{Log}
    \DeclareMathOperator{\Exp}{Exp}

	\DeclareMathOperator{\arccosh}{arccosh}

    \newcommand{\half}{\frac{1}{2}}

    \newcommand{\define}[1]{\emph{#1}}      

    \newcommand{\ltwo}{l^{2}}
    \newcommand{\Ltwo}{L^{2}}
    
    \newcommand{\VN}[1]{\mathcal{#1}}

    \newcommand{\len}{\ell}   
    
	\newcommand{\Del}{\Delta}
	
\begin{document}

%
%
\title{Zeta functions of graphs with $\Z$ actions}

\author{Bryan Clair}
\address{Saint Louis University,
         220 N. Grand Avenue
		 St. Louis, MO 63108}
\email{bryan@slu.edu}
\date{\today}


\begin{abstract}
Suppose 
$Y$ is a regular covering of a graph $X$ with covering transformation 
group $\pi = \mathbb{Z}$.
This paper gives an explicit formula for the $L^2$
zeta function of $Y$ and computes examples.
When $\pi = \mathbb{Z}$, the $L^2$ zeta function is
an algebraic function.  As a consequence it extends to
a meromorphic function on a Riemann surface.
The meromorphic extension provides a setting to generalize
known properties of zeta functions of regular graphs, such as
the location of singularities and the functional equation.
\end{abstract}

\maketitle

\section{Introduction}
Given a finite graph, there is a zeta function which encodes some of 
the combinatorics of the graph.  The zeta function was defined by Ihara and 
extended by Hashimoto and then Bass.  See~\cite{starkter, hst:zeta} 
for a fine introduction to the subject taking a geometric approach.

There is an analogous zeta function for any infinite graph with 
cofinite action of a discrete group. 
Let $Y = (VY,EY)$ be a locally finite (but typically infinite) graph and
suppose the group $\pi$ acts freely on $Y$ with finite quotient graph $X$.
Let $P$ denote the set of
free homotopy classes of primitive closed paths in $Y$.
For $\gamma\in P$, $\len(\gamma)$ is the length of the shortest
representative of $\gamma$.
The group $\pi_{\gamma}$ is the stabilizer of $\gamma$ under the
action of $\pi$.
The \define{$\Ltwo$ zeta function} of $Y$ is the infinite product
\begin{equation}\label{eq:l2z}
   Z^{(2)}_{Y}(u)^{-1} = \prod_{\gamma \in \pi\backslash P}
				   \left(1-u^{\len(\gamma)}\right)
				   ^{\frac{1}{\abs{\pi_{\gamma}}}}.
\end{equation}
  This definition was first given in~\cite{cms:zeta} as a 
specialization from a more general setting, but beware that 
the notation $Z$ in~\cite{cms:zeta} refers to the reciprocal of the zeta
function considered here and elsewhere in the literature.  
See~\cite{gil:zeta} for a more direct treatment of the case 
considered here.

  For finite graphs, the fundamental theorem is the Ihara-Hashimoto-Bass 
rationality formula, which says that the zeta function is the 
reciprocal of a polynomial.  The analogous theorem for infinte graphs 
requires techniques of von Neumann algebras.  The infinte graph result
is formally similar, and implies
convergence of \eqref{eq:l2z}, but the $\Ltwo$ zeta function is not 
typically a rational function.

Let $\delta$ be the adjacency operator of $Y$ acting on
$\ltwo (VY)$. For $f\in \ltwo (VY)$ let $Qf(v)=q(v)f(v)$ where $q(v)+1$
is the degree of the vertex $v$.
Put $\Del_{u}=I-\delta u + Q u^2$.
Then from \cite[Theorem 0.3]{cms:zeta},
\begin{equation}\label{eq:zetadef}
   Z_{Y}^{(2)}(u)^{-1} = (1-u^{2})^{-\chi(X)}\Det_{\pi}\Del_{u}
\end{equation}
where $\Det_{\pi}$ is a von Neumann determinant defined 
in~\cite{cms:zeta}.
In particular, the product in \eqref{eq:l2z} converges for small $u$
(which was not a priori obvious).  

In this paper, the only group considered is $\pi = \Z$, so that
$X = Y/\Z$.
Theorem~\ref{thm:z} computes $Z_{Y}^{(2)}(u)$ in this case.  The main 
difficulties to overcome are the evaluation of a particular definite 
integral and careful bookeeping with branches of multi-valued complex 
functions.

The formula for $Z_{Y}^{(2)}(u)$ is algebraic, and 
Theorem~\ref{thm:extend}
takes advantage of this to extend $Z_{Y}^{(2)}(u)$ to a meromorphic 
function $\tilde{Z}$ defined on a compact
Riemann surface $S$ (which depends on $Y$).
From another 
viewpoint, $Z_{Y}^{(2)}(u)$ is naturally a multi-valued meromorphic
function defined on all of $\C$.

The  surface $S$ covers the Riemann sphere $\C P^{1}$ with 
branch points, and the branch points play a similar role for
infinite graphs as the poles do for zeta functions of finite graphs.
Specifically, Theorem~\ref{thm:branchpts} gives conditions for 
$\tilde{Z}$ of a 
$q+1$ regular graph $Y$ to have all its branch points over the set
\[
   C=\{u\in \C : \abs u=q^{-1/2} \} \cup [ -1 , - \frac{1}{q}] 
	   \cup [\frac{1}{q} , 1].
\]
$C$ is  exactly the set where poles may occur for zeta 
functions of finite $q+1$ regular graphs.

The extension to $\tilde{Z}$ gives a meaningful context for 
functional equations relating $u  \leftrightarrow \frac{1}{qu}$, and 
Section~\ref{sec:func} explores these.

Finally, Section~\ref{sec:examples} gives a number of computations 
for specific $Y$.

This paper is intended as a model for how one might attack more 
general $\pi \neq \Z$.
It is shown in~\cite{cms:zeta2} that
for a $q$-regular graph, the $\Ltwo$ zeta function always extends
holomorphically to the interior of the set $C$.
In the most optimistic scenario, the $\Ltwo$ 
zeta function is always algebraic and therefore extends past $C$
to a compact Riemann surface.
More likely, one may need to allow noncompact surfaces with 
infinitely many sheets over $\C P^{1}$.  In the worst scenario,
the ``branch points'' could spread out continuously over $C$
and prevent any further extension of domain.
In any event, the explicit computation that provides the key here
is not likely to unlock the more general case.

The author would like to thank David Moulton and Shahriar Mokhtari
for helpful discussions.

\subsection{Group Von Neumann Algebras}
For completeness, here is a quick overview of relevant material from  
von Neumann algebras.
For $\pi$ a countable discrete group, the \define{von Neumann algebra}
of $\pi$ is the algebra $\VN{N}(\pi)$ of bounded $\pi$-equivariant operators
from $\ltwo(\pi)$ to $\ltwo(\pi)$.

The \define{von Neumann trace} of
an element $f \in \VN{N}(\pi)$ is defined by
\[ \Tr_{\pi}f = \inner{f(e)}{e} \]
for $e \in \pi$ the unit element.
The group ring $\C[\pi]$ is contained in $\VN{N}(\pi)$, acting on 
$\ltwo(\pi)$ by right multiplication.  It is a dense subspace.
The trace of an element of the 
group ring is simply the coefficient of the identity.

For $H =
\oplus_{i=1}^{n}\ltwo(\pi)$ and a bounded $\pi$-equivariant
operator $f:H \to H$,
define
\[ \Tr_{\pi}f = \sum_{i=1}^{n}\Tr_{\pi}f_{ii}. \]
The trace as defined is independent of the decomposition of $H$.
The determinant $\Det_{\pi}\Del(Y,u)$ is defined via formal power
series as $(\Exp \circ \Tr_{\pi} \circ \Log) \Del(Y,u)$ and 
converges for small $u$.

\begin{Example}
	When $\pi$ = $\Z$, Fourier transform identifies $\ltwo(\Z)$ with 
	$\Ltwo(S^{1})$.  An element 
	$\sum_{n=-\infty}^{\infty}c_{n}t^{n} \in \VN{N}(\pi)$ tranforms to 
	multiplcation by
	$f(\theta) = \sum_{n=-\infty}^{\infty}c_{n}e^{in\theta}$, and
	\[
	\Tr_{\pi}f = \inner{f\cdot 1}{1} = \int_{S^{1}}f(\theta)d\theta = 
	c_{0}.
	\]
\end{Example}

\section{Graphs with $\Z$ actions}
We assume $\pi = \Z = \left<t\right>$, the free abelian group on one generator 
$t$.  Suppose $X = Y/\Z$ has $v$ 
vertices.  Choosing lifts of these vertices to $Y$, we identify 
\[
  \ltwo(VY) = \bigoplus_{v}\ltwo(\Z)
\]
and the adjacency operator $\delta$ is then a $v \times v$ matrix with
entries in the group ring $\Z[\left<t\right>]$.
Since $\delta$ is self-adjoint, it satisfies $\delta(t) = 
\delta(t^{-1})^{T}$.  Similarly, $\Del_{u}(t) = 
\Del_{u}(t^{-1})^{T}$ (but beware that $\Del_{u}$ is not generally 
self-adjoint).  Therefore, $\Det \Del_{u} \in \C[\left<t\right>]$ is 
symmetric in $t$ and $t^{-1}$, and we can write $\Det \Del_{u}  = 
P_{u}(\frac{t+t^{-1}}{2})$ for some polynomial $P_{u}$.  The 
coefficients of $P_{u}$ are integer polynomials in $u$.

We know in general that $\Det_{\pi}  \Del_{u}$ is idependent of the
choice of lifts of vertices, but here 
it is very clear, since choosing a different lift
will multiply a row by $t^{k}$ and the corresponding column by 
$t^{-k}$ (for some $k$).  In particular, $P_{u}$ depends only
on $Y$ and the $\Z$ action.

Now, under Fourier transform, $\bigoplus_{v}\ltwo(Z) = 
\bigoplus_{v}\Ltwo(S^{1})$.  Here $S^{1} = \{e^{i\theta} | \theta 
\in (-\pi,\pi]\}$ with measure normalized to have total measure 1.
Under Fourier tranform, multiplication by $t$ becomes
multiplication by the function $e^{i\theta}$, and hence $\Del_{u}$ is 
represented by a $v \times v$ matrix which will be denoted 
$M_{u}(\theta)$.
To compute the zeta function,
\begin{align}
	\Det_{\pi} \Del_{u} &= \exp \Tr_{\pi} \Log(\Del_{u}) \\
	 &= \exp \int_{S^{1}} \Tr \Log(M_{u}(\theta)) d\theta \\
	&= \exp \int_{S^{1}} \log \det(M_{u}(\theta)) d\theta \\
	&= \exp \int_{S^{1}} \log P_{u}(\cos(\theta)) d\theta. 
	\label{eq:detlog}
\end{align}

\subsection{The Line}\label{sec:line}
To proceed further, we work out a crucial example.
Let $VY = \Z$, and connnect $n$ to $n+1$ with an edge (so $Y$ is a 
line).  Then
\[ \Del_{u} = 1 - (t + t^{-1})u + u^{2},\]
\[ M_{u}(\theta) = 1 - 2\cos(\theta)u + u^{2}, \]
and
\[ P_{u}(x) = 1 - 2ux + u^{2}.\]
Notice that for $\abs{u} < 1$ and for all $\theta$, $M_{u}(\theta) 
\notin (-\infty,0]$. In what follows, $\log$ will be the principal 
branch of the logarithm.

Now restrict to $\abs{u} < 1$.  Because $Y$ has no loops,
the $\Ltwo$ zeta function for $Y$ is 
identically 1.  Therefore, by \cite[Theorem 0.3]{cms:zeta}
\begin{align}
	1 &= (1-u^{2})^{0}\Det_{\pi}\Del_{u}
	\nonumber \\
	&= \exp \int_{S^{1}}\log (1 - 2u\cos(\theta) + u^{2}) d\theta \\
	&= \exp \int_{S^{1}}\log(2u) + \log(\frac{u+u^{-1}}{2} - 
	   \cos(\theta)) d\theta \\
	&= 2u \exp \int_{S^{1}} \log(r - \cos(\theta)) d\theta.
\end{align}
Here, we have assumed $u \neq 0$ and put $r = (u+u^{-1})/2$.  
Generally, some 
care must be taken when writing $\log(xy) = \log(x) + \log(y)$.  
If $u \in (-1,0)$, then $r - \cos(\theta) < 0$ and the identity
is off by $2\pi i$.  However, 
the $2\pi i$ is washed out by the $\exp$ in front.  For other values of $u$
there is no problem, because the imaginary parts of 
$u$ and $r$ have opposite sign.  

Notice that $r = \cosh (-\log(u))$, so that $u = e^{-\arccosh(r)}$.  
Here, $\arccosh$ has a branch cut discontinuity on $(-\infty,1]$ and 
range $\{a+bi|a>0, b \in (-\pi,\pi]\} \cup [0,\pi]i$.

\begin{Prop}\label{prop:integral}
	Suppose $r \in \C$.  Then
	\begin{equation}\label{eq:int}
		\int_{S^{1}} \log(r - \cos(\theta)) d\theta =
		   \arccosh(r) - \log(2).
    \end{equation}
\end{Prop}
\begin{proof}
The discussion above proves that for $r \in \C-[-1,1]$,
\[
  \exp \int_{S^{1}} \log(r - \cos(\theta)) d\theta = \half e^{\arccosh(r)}.
\]
Taking the $\log$ of both sides,
\begin{equation}\label{eq:notquite}
  \int_{S^{1}} \log(r - \cos(\theta)) d\theta =
  \arccosh(r) - \log(2)
\end{equation}
In principle, this is only true up to $2\pi i k$ for some $k \in \Z$.
However, $k$ must be zero 
since both $\log$ and $\arccosh$ have imaginary part in the range 
$(-\pi,\pi]$.

Now we extend \eqref{eq:notquite} to all of $\C$.  We check the 
imaginary part explicitly.  For $r \in [-1,1]$, put $r = 
\cos(\phi)$, $\phi \in [0,\pi]$.  Then $\Im(\arccosh(r) - \log 2) = 
\phi$.  On the other hand, $\arg(r-\cos(\theta)) = \pi$ when 
$\cos(\theta) > r$ and is 0 otherwise.  Therefore,
\begin{align}
  \int_{S^{1}}\Im(\log(r-\cos(\theta))) d\theta &=
  \int_{S^{1}}\arg(r-\cos(\theta)) d\theta\\
  &=  \pi \cdot m\{\theta | \cos(\theta) > \cos(\phi)\}\\
  &=  \pi \cdot \frac{2\phi}{2\pi} =  \phi.
\end{align}

Next, consider the real part of \eqref{eq:notquite}.
The real part of the left hand side is
\[
   \int_{S^{1}}\log\abs{r-\cos(\theta)} d\theta.
\]
It is not hard to see that this integral is finite even for
$r \in [-1,1]$.  On the other hand, $\Re(\arccosh(r) - \log 2)$ is
a continuous function on all of $\C$  (and equals $-\log 2$ on 
$[-1,1]$).  Thus the two sides are
defined on all of $\C$, equal on $\C - [-1,1]$, and the right side
is continuous.

Now, for $r \in [-1,1], n=1,2,\ldots$ put
\[
  f_{n}(\theta) = \log \abs{r - \cos(\theta) + i/n}.
\]
The $f_{n}$ are a decreasing sequence of functions, bounded above 
(by $\sqrt{5}$), and converging a.e. to $\log\abs{r-\cos(\theta)}$.  By 
the Monotone Convergence Theorem,
\begin{align}
  -\log 2 &= \lim_{n \to \infty}\Re(\arccosh(r + i/n) - \log 2)\\
          &= \lim_{n \to \infty} \int_{S^{1}}f_{n}(\theta)d\theta\\
          &=\int_{S^{1}}\log\abs{r-\cos(\theta)} d\theta,
\end{align}
\end{proof}

\begin{Remark}
	Computing the integral in \eqref{eq:int} is a good, difficult 
	calculus exercise for $r > 1$.  I know of no elementary way to compute 
	it in general.
\end{Remark}

\begin{Remark}
	The inverse hyperbolic cosine function satisfies
\[  \arccosh(r) = \log \left(r + \sqrt{r+1}\sqrt{r-1}\right) \]
where the principal branches of $\arccosh$, $\log$, and $\sqrt{z}$
are used.  In particular,
$\Re(\arccosh(r)) = \log \Abs{r + \sqrt{r+1}\sqrt{r-1}}$ is a 
continuous function.
Taking the real part of both sides of \eqref{eq:int} gives
the integral
	\begin{equation}\label{eq:absint}
		\int_{S^{1}}\log\abs{r-\cos(\theta)} d\theta =
		\log \half\Abs{r + \sqrt{r+1}\sqrt{r-1}}.
	\end{equation}
for all $r \in \C$.
\end{Remark}

\subsection{The explicit formula}
\begin{Thm}\label{thm:z}
	Let $Y$ be a regular $\Z$ covering of a finite graph $X$.  Let $P_{u}(x)$ 
	be the degree $n$ polynomial so that
	\[ \det \Del_{u} = P_{u}\left(\frac{t+t^{-1}}{2}\right). \]
	There is $R > 0$ so that for all $0 < \abs{u} < R$,
	\begin{equation}
		Z_{Y}^{(2)}(u)^{-1} = 
		(1-u^{2})^{-\chi(X)}\frac{\alpha(u)}{2^{n}}
		\prod_{i}\left(r_{i} + \sqrt{r_{i}+1}\sqrt{r_{i}-1}\right).
	\end{equation}
	Here $P_{u}(x) = \alpha(u)\prod_{i=1}^{n}(r_{i}(u)-x)$, 
	and $r_{i}(u)$ are the roots of $P_{u}$.  The square roots are 
	principal, in the sense that $\sqrt{z} = \exp(\half\log(z))$.
\end{Thm}
\begin{proof}
	The polynomial $(-1)^{n}\alpha(u)$ is the coeffieicent of the top degree 
	term $x^{n}$ of $P_{u}$.  Since $P_{0} = 1$, 0 is a root of 
	$\alpha$.  There is $R_{1} > 0$ with $\alpha(u)\neq 0$ on 
	$0 < \abs{u} < R_{1}$, and so one can write
	$P_{u}(x) = \alpha(u)\prod_{i=1}^{n}(r_{i}(u)-x)$.
	
	From~\eqref{eq:zetadef}, one need only compute $\Det_{\pi} \Del_{u}$.
	There is a subtle point involving the log of a product, but the 
	heart of the argument is the computation below,
	which begins with~\eqref{eq:detlog}, and uses Proposition~\ref{prop:integral}:
\begin{align}
	\Det_{\pi} \Del_{u}
	&= \exp \int_{S^{1}} \log P_{u}(\cos(\theta)) d\theta \\
	&= \exp \int_{S^{1}} \log 
       \alpha(u)\prod_{i=1}^{n}(r_{i}(u)-\cos(\theta))d\theta 
	   \label{eq:log1} \\
    &= \exp \left( \log \alpha(u) +
	               \sum_{i=1}^{n}\int_{S^{1}}
				   \log (r_{i}(u)-\cos(\theta)) d\theta \right)
       \label{eq:log2} \\
	&= \exp \left( \log \alpha(u) +
	               \sum_{i=1}^{n}(\arccosh(r_{i})-\log(2)) \right)\\
	&= \frac{\alpha(u)}{2^{n}}\prod_{i}\exp(\arccosh(r_{i}))\\
	&= \frac{\alpha(u)}{2^{n}}
		\prod_{i}\left(r_{i} + \sqrt{r_{i}+1}\sqrt{r_{i}-1}\right).
\end{align}

    It remains to justify the transition from~\eqref{eq:log1} 
	to~\eqref{eq:log2}.

	Write
	\begin{equation}\label{eq:kdef}
		\log P_{u}(x) =  \log \alpha(u) +
	                      \sum_{i=1}^{n}\log (r_{i}(u)-x) +
						  2\pi i k(u,x).
	\end{equation}
	The function $k(u,x)$ is always an integer.  We will show that
	$k(u,x) = k(u)$ is independent of $x \in [-1,1]$ and therefore 
	pulls through the integral in \eqref{eq:log1} to be eaten by the 
	$\exp$.
	
	Since $\Delta_{u} = I-\delta u+Q u^{2}$, we can write $P_{u}(x) = 
	1 + u T_{u}(x)$ for some polynomial $T$.
	Then there is $R_{2}>0$ so that for $\abs{u} < 
	R_{2}$ and $x \in [-1,1]$ we have $\Re(P_{u}(x)) > 0$. 
	Therefore, $\log(P_{u}(x))$ is a continuous function of $x \in 
	[-1,1]$.
	
	In addition, for $0 < \abs{u} < R_{2}$, we see that $P_{u}(x)$ 
	has no roots on $[-1,1]$, \ie\ $r_{i}(u) \notin [-1,1]$.
	Therefore $\log(r_{i}(u) - x)$ is a continuous function of $x \in 
	[-1,1]$ (since we're using the principal branch of the logarithm).
	
	We have shown that all other terms in~\eqref{eq:kdef} are continuous 
	functions of $x$, and therefore $k(u,x)$ is a continuous function 
	of $x$ on $[-1,1]$, hence constant in $x$.
	
	Setting $R = \min\{R_{1},R_{2}\}$ completes the proof.
\end{proof}

\subsection{The meromorphic extension}
From Theorem~\ref{thm:z}, it is apparent that $Z_{Y}^{(2)}(u)$ is an 
algebraic function of $u$.  In this section, we make this more 
explicit and then explore the consequences.

Let $s_{i}= \sqrt{r_{i}+1}\sqrt{r_{i}-1}$, and for 
$I=(\iota_{1},\ldots,\iota_{n}) \in \{\pm 1\}^{n} = \Z_{2}^{n}$, put
\[
   W_{I} = \prod_{i=1}^{n}r_{i}+\iota_{i}s_{i}.
\]
Note that
$(r_{i}+s_{i})(r_{i}-s_{i})=1$ so that $W_{I}^{-1}=W_{-I}$.
Theorem~\ref{thm:z} then says that
\begin{equation}\label{WfromZ}
	Z_{Y}^{(2)}(u) = (1-u^{2})^{\chi(X)}\frac{2^{n}}{\alpha(u)}
	   W_{-1,-1\ldots ,-1}.
\end{equation}
Let
\begin{equation}\label{eq:omegadef}
  \Omega(T) = \prod_{I \in \Z_{2}^{n}}(T-W_{I}).
\end{equation}
Then $\Omega$ is a polynomial in $T$ of degree $2^{n}$.  It is 
invariant under the transformation $s_{i}\to -s_{i}$, hence it is 
even degree in each $s_{i}$.  We can replace $s_{i}^{2}$ with 
$r_{i}^{2}-1$ so that $\Omega$ is a degree $n$ polynomial in $r_{i}$, symmetric 
in the $r_{i}$.  This means that $\Omega$ is in fact a polynomial in
the elementary symmetric functions $\sigma_{1},\ldots,\sigma_{n}$ of 
the $r_{i}$, for example:
\begin{equation}\label{eq:omega}
   \Omega(T)=
   \begin{cases}
	   T^2-  2\,T\,{{\sigma }_1}  + 1
	      \quad & \text{for $n = 1$,}\\
	   T^4-  4\,T^3\,{{\sigma }_2}  +
	       T^2\,\left( -2 + 4\,{{{\sigma }_1}}^2 -
		  8\,{{\sigma }_2} \right)   -
		   4\,T\,{{\sigma }_2}  + 1
		  \quad & \text{for $n = 2$,}
    \end{cases}
\end{equation}
and when $n=3$,
\begin{equation}
\begin{split}
   \Omega(T)=
	  T^8 &-  8T^7{\sigma_3}  + 
		 T^6\left( 4 - 8{{\sigma_1}}^2 + 
			16{\sigma_2} + 16{{\sigma_2}}^2 - 
			32{\sigma_1}{\sigma_3} \right)  \\
		 &-  T^5\left( -40{\sigma_3} + 
			32{{\sigma_1}}^2{\sigma_3} - 
			64{\sigma_2}{\sigma_3} \right)   \\
		 &+  T^4\left( 6 - 16{{\sigma_1}}^2 + 
			16{{\sigma_1}}^4 + 32{\sigma_2} - 
			64{{\sigma_1}}^2{\sigma_2} + 
			32{{\sigma_2}}^2 + 
			64{\sigma_1}{\sigma_3} + 
			64{{\sigma_3}}^2 \right)\\
		 &- \cdots - 8T{\sigma_3} + 1
\end{split}
\end{equation}
using the symmetry of coefficients to finish (roots
of $\Omega$ occur in reciprocal pairs).

Since the $r_{i}$ are the roots of $P_{u}$,
\[
  \sigma_{i} = (-1)^{n-i}\left(
    \text{the $n-i^{\text{th}}$ coefficient of\ }
	\frac{P_{u}}{\alpha(u)} \right).
\]
Thus $\sigma_{i}$ is a rational
function of $u$, and so $\Omega \in \C(u)[T]$.

We have shown that $W_{I}$ and therefore $Z_{Y}^{(2)}(u)$
are algebraic functions of $u$ of degree less than or equal to 
$2^{n}$.

\begin{Thm}\label{thm:extend}
Let $Y$ be a regular $\Z = \pi$
covering of a finite graph $X$.
Then $Z_{Y}^{(2)}(u)$ extends uniquely
to a meromorphic function on a Riemann surface.

More precisely, there exists a compact Riemann surface $S$, a (branched) 
covering map $\Pi : S \to \C P^{1}$, and a meromorphic function
$\tilde{Z}$ on $S$.  There is a point $z_{0} \in \Pi^{-1}(0)$
and a neighborhood $U$ of $z_{0}$ on which $\Pi$ is biholomorphic
such that $\tilde{Z}(z) = Z_{Y}^{(2)}(\Pi(z))$ for all $z \in U$.

The triple $(S,\Pi,\tilde{Z})$ is unique in the following sense:
If $(S^{\prime},\Pi^{\prime},\tilde{Z}^{\prime})$ has the 
corresponding properties, then there exists exactly one
fiber preserving biholomorphic mapping $\tau : S \to S^{\prime}$
such that $\tilde{Z} = \tilde{Z}^{\prime}\circ \tau$.
\end{Thm}
\begin{Remark}
	The number of sheets of $\Pi$ is less than or equal to $2^{n}$, where 
	$n$ is the degree of the polynomial $P_{u}$ defined earlier.
\end{Remark}

\begin{proof}
Define $W_{I}$ and $\Omega$ as above.
The difficult work is finished, as we showed already that 
$W_{I}$ is algebraic.
Since $W_{-1,-1,\ldots,-1}$ is holomorphic in a neighborhood of 0
and $\Omega(W_{-1,-1,\ldots,-1}) = 0$, there is a unique irreducible
factor $\Phi \in \C(u)[T]$ 
with 
\begin{equation}\label{eq:Phi}
	\Phi(W_{-1,-1,\ldots,-1})=0.
\end{equation}
in a neighborhood of 0.

The algebraic function defined by $\Phi(T)$ consists of $S$ and $\Pi$ 
as above, plus a meromorphic function $f$ on $S$ such that
$(\Pi^{*}\Phi)(f) = 0$.  
It is unique in the sense of fiber preserving biholomorphic mappings
as above (see~\cite[I.8]{forster} for details).

Since $W_{-1,-1,\ldots,-1}$ is holomorphic 
in a neighborhood of 0, there is a point $z_{0} \in \Pi^{-1}(0)$
and a neighborhood $U$ of $z_{0}$ on which $\Pi$ is biholomorphic
with $f(z) = W_{-1,-1,\ldots,-1}(\Pi(z))$ for $z \in U$.

For $z \in S$, let $u = \Pi(z)$ and put
\begin{equation}\label{eq:ztilde}
      \tilde{Z}(z) = (1-u^{2})^{\chi(X)}\frac{2^{n}}{\alpha(u)}f(z).
\end{equation}
to complete the proof.
\end{proof}

\section{Regular graphs}
In this section, assume that $X$ is $q+1$ regular.

\subsection{Functional Equations}\label{sec:func}
The zeta function for finite regular graphs satisifes a number of
functional equations under the transformation
\[
   \tau : u \to \frac{1}{qu}
\]
(see~\cite{starkter}).  The situation for $\Ltwo$ zeta functions is 
somewhat less simple.

First notice that
\[
   \Del_{1/qu} = I - \delta\frac{1}{qu} + q\frac{1}{(qu)^{2}}
               = \frac{1}{qu^{2}}\left(I - \delta u + qu^{2}\right)
			   = \frac{1}{qu^{2}}\Del_{u}.
\]
Then the polynomial $P_{1/qu}(x)$ has the same roots 
$r_{1},\ldots,r_{n}$ as $P_{u}(x)$.
Since $\Omega$ and $W_{-1,-1,\ldots,-1}$ are symmetric functions of 
the $r$'s, they are invariant under $\tau$.

Suppose $\Omega$ is irreducible, so that the $\Ltwo$ zeta function is 
defined on the Riemann surface $S$ for $\Omega$ by~\eqref{eq:ztilde}.
Then the transformation $u \to \frac{1}{qu}$
induces a biholomorphic involution $\tilde{\tau}:S \to S$ so that
$f \circ \tilde{\tau} =  f$.  It is then easy to find functional
equations for $\tilde{Z}$.  For example:
\begin{Prop}
	Suppose $X$ is $q+1$-regular and $\Omega$ is irreducible.
	For $z \in S$ put $u = \Pi(z)$.
	Then
	\begin{equation}\label{eq:func}
		\left(\tilde{Z}\circ \tilde{\tau}\right)(z) =
		  q^{2e-v}u^{2e}{\left(\frac{1-u^{2}}{q^{2}u^{2}-1}\right)}^{-\chi}
		  \tilde{Z}(z).
	\end{equation}
	Here $v$ and $e$ are the number of vertices and edges of $X$, and
	$\chi = \chi(X) = v-e$.
	(Compare~\cite[Cor 3.10]{bass:zeta})
\end{Prop}
\begin{proof}
	This is a straightforward calculation using
	$f \circ \tilde{\tau} =  f$, equation~\eqref{eq:ztilde}, and
	\[ \alpha(\frac{1}{qu}) = 
	{\left(\frac{1}{qu^{2}}\right)}^{v}\alpha(u). \]
\end{proof}

If $\Omega$ is reducible, one gets a collection of disjoint Riemann 
surfaces $S_{1},\ldots,S_{k}$ and the map $\tilde{\tau}$ may permute
them.  We are interested in $\tilde{Z}$ on a particular choice $S$,
and so it will not satisfy a functional equation in any traditional 
sense.  The line (example~\ref{ex:line}) is a good example of this.

\subsection{Location of branch points}
The zeta function for a finite, $q+1$ regular graph has all of its 
poles in the set
\[
   C=\{u\in \C : \abs u=q^{-1/2} \} \cup [ -1 , - \frac{1}{q}] 
	   \cup [\frac{1}{q} , 1].
\]
For the $\Ltwo$ zeta function, we can make a slightly weaker 
statement for branch points.

\begin{Thm}\label{thm:branchpts}
	Let $Y$ be a regular $\Z = \pi$ covering of a finite graph $X$.
	Suppose that $X$ is $q+1$ regular.  Let $\Omega$ be the polynomial 
	defined in~\eqref{eq:omegadef}, and assume $\Omega$ is irreducible.
	If the field extension
	$\C(u)[T]/(\Omega(T)) : \C(u)$ is Galois, then the covering
	$\Pi$ from Theorem~\ref{thm:extend} has all of its branch points
	over $C$.
\end{Thm}
\begin{proof}
	Let $D_{0}$ and $D_{\infty}$ be the connected components of $\C - C$.
	From~\cite{cms:zeta2}, the
	$\Ltwo$ zeta function $Z_{Y}^{(2)}(u)$ extends holomorphically
	to $D_{0}$, so the neighborhood $U$ from Theorem~\ref{thm:extend}
	must also extend to cover $D_{0}$ with no branch points.
	The field extension
	$\C(u)[T]/(\Phi(T)) : \C(u)$ is Galois if and only if the deck
	transformations of $S$ over $\C P^{1}$ act transitively on the
	sheets of $S$ (\cite[pg. 57]{forster}).  Then $\Pi^{-1}(D_{0})$ is
	a union of copies of $U$ and has no branch points.
	The involution $\tilde{\tau}$ from the functional equation
	biholomorphically interchanges 
    $\Pi^{-1}D_{0}$ with $\Pi^{-1}D_{\infty}$,
	so that $\Pi$ can only be branched on $C$.
\end{proof}

In example~\ref{ex:sawtooth}, we will see a graph for which the
$\Ltwo$ zeta function is branched over $0$ and the deck 
transformations of $S$ are not transitive.

The assumption that $\Omega$ is irreducible is less well motivated.
As in example~\ref{ex:triladder}, the zeta function for a graph with
reducible $\Omega$
will still satisfy a functional equation if $\tilde{\tau}$ perserves $S$.

The following argument gives hope for a close relationship between
branch points of $\tilde{Z}$ for $Y$ and zeros of $Z(X)$.
To compute the zeta 
function $Z(X)$ of the quotient graph $X = Y/\Z$, one takes the
determinant of $\Delta_{X}(u) = I -\delta_{X}u + Q u^{2}$, where
$\delta_{X}$ is the adjacency operator on $X$.
Poles of $Z(X)$ occur when $\det \Delta_{X}(u) = 0$.
But $\delta_{X}$ is equal to $\delta$ on $Y$ under $t\to 1$,
and so poles of $Z(X)$ occur when $P_{u}(1)=0$, or equivalently
when some root $r_{i}(u) = 1$.

If $r_{i}(u) = 1$ then the terms
$r_{i} \pm \sqrt{r_{i}+1}\sqrt{r_{i}-1}$ coincide.  In other words,
two roots of $\Omega$ coincide at any $u$ where $Z(X)$ has a pole --
a necessary condition for $S$ to be branched over $u$.

Frequently, branch points of $\tilde{Z}$ do coincide with poles
of $Z(X)$.  However, examples in the next section show that
both possible implications are false in general.

\section{Examples}\label{sec:examples}
\begin{Example}[The Line]\label{ex:line}
Let $Y$ be the line, as in Section~\ref{sec:line}.  We saw earlier 
that
\[ P_{u}(x) = 1+u^{2}- 2ux.\]
Then $\alpha(u) = 2u$ and $r(u)=\frac{1+u^{2}}{2u}$.
From~\eqref{eq:omega}, we have
\[
   \Omega(T) = T^{2} - T\frac{1+u^{2}}{u} + 1
   = \frac{(T-u)(Tu-1)}{u}.
\]
Here $\Omega$ is reducible. Some careful computation shows that
\begin{align}
   W_{-1}(u) &= \frac{1+u^{2}}{2u} - 
   \sqrt{\frac{1+u^{2}}{2u}+1}
   \sqrt{\frac{1+u^{2}}{2u}-1}
   &=\begin{cases}
         u \quad \text{if $\abs{u} < 1$}\\
		 1/u \quad \text {if $\abs{u} > 1$}
	 \end{cases}
\end{align}
so $\Phi(T) = T-u$, the Riemann surface $S$ is $\C P^{1}$, $f(u)=u$,
and the zeta function is $2f/\alpha = 1$.

Notice that the transformation $\tau : u \to \frac{1}{qu}$
(here $q=1$) interchanges the two irreducible surfaces.  On the other 
surface, the analog of $\tilde{Z}$ is $u^{2}$, and in fact the 
functional equation~\eqref{eq:func} becomes
\[
   u^{2} = 1^{2\cdot 1 - 1}u^{2\cdot 1}\cdot 1 \cdot 1.
\]
\end{Example}
\begin{Example}[Some degree 1 graphs]\label{ex:deg1}
Let $Y$ be the first graph shown in Table~\ref{deg1table} (all these 
graphs take the obvious $\Z$ action).  $Y$ is 
4-regular, so $q=3$.  It's quotient graph $X$ is a vertex with two
loops,and $\chi(X) = -1$.

The adjacency matrix for $Y$ is the $1\times 1$ matrix
$(t^{-1}+2+t)$.  Then $P(x) = -2ux + 1 - 2u + 3u^{2}$ which has the 
one root shown in the table.  From~\eqref{eq:omega},
\[ \Omega(T) = T^{2} - \left(\frac{1-2u+3u^{2}}{u}\right)T + 1 \]
which is irreducible.
 
The associated Riemann surface $S$ is a two sheeted branched cover
of $\C P^{1}$.  Possible branch points occur when the 
discriminant of $\Omega$ vanishes, which happens in this case at
\[ u=1, u=\frac{1}{3}, u=\frac{i}{\sqrt{3}}, u=\frac{-i}{\sqrt{3}}. \]
Here, all four are in fact branch points of multiplicity 2.
The pattern of branch points is shown the table, and the set $C$ is 
also indicated.

The Riemann-Hurwitz formula gives the genus of a 
branched covering of $\C P^{1}$ as $g = b/2 - d + 1$ with $d$ the 
number of sheets and $b$ the total branching order.
For this graph the genus is 1 and $S$ is a torus.

Other lines of Table~\ref{deg1table} give the results of similar 
computations for different $Y$ with $n=1$.
In all cases, $\Omega(T) = T^{2} - 2rT + 1$ is 
irreducible, $S$ is a two sheeted branch cover, and all branch points 
are mulitiplicity 2.

Graph \#3 is an example in which poles of the zeta function for the 
quotient graph do not correspond to branch points of $S$.
In this graph \#3 of the table, $r(- \frac{1}{4} \pm 
\frac{i}{4}\sqrt{7}) = 1$,
but these are not branch points of the $\Ltwo$ zeta function.

Graphs \#2,4, and 5 are bipartite and have vertical bilateral symmetry.
Graphs \#4 and 5 have different zeta functions because they have
different $\alpha$ and different $\chi$.

Graph \#6 is non-regular.  It's branch points are shown with
circles of radius $1/\sqrt{2}$ and $1/\sqrt{3}$ for scale.
\end{Example}

\begin{table}\label{deg1table}
\begin{tabular}{|l|c|c|c|c|c|c|}
	\# & Graph & $q$ & $\chi(X)$ & $r(u)$ & Branchpoints & Genus \\
	\hline
	1 &	\includegraphics[width=1in]{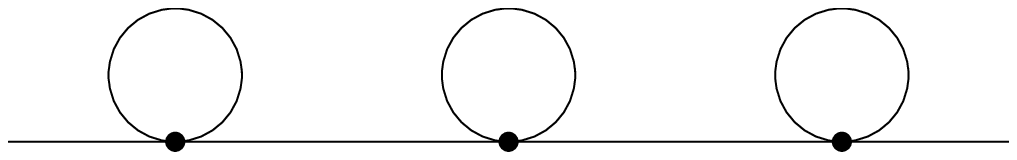} 
	   & 3 & -1
	   & $\frac{1 - 2\,u + 3\,u^2}{2u}$
	   & 	\includegraphics[width=1in]{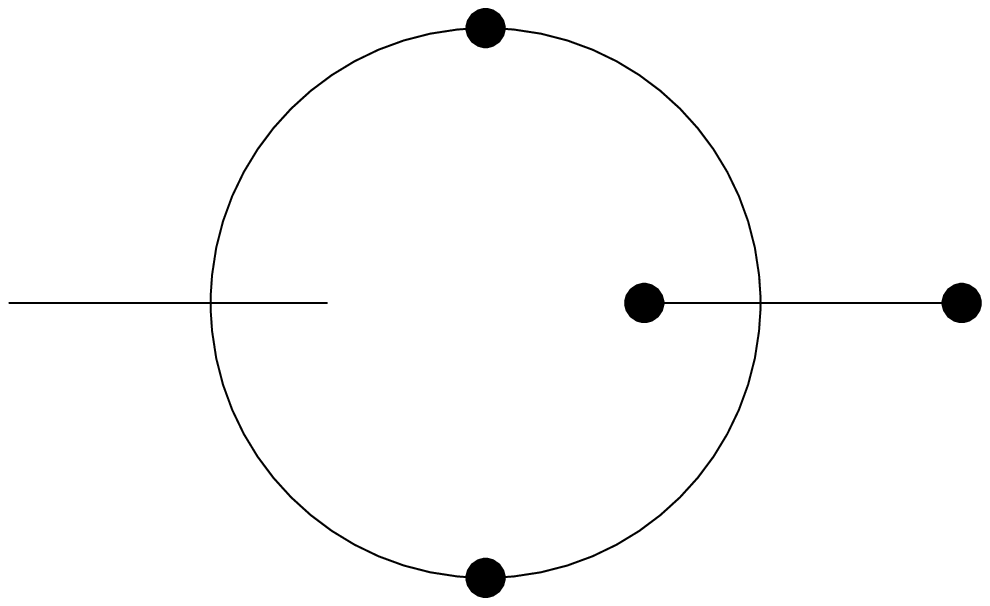} 
	   & 1 \\
	   \hline
	2 & \includegraphics[width=1in]{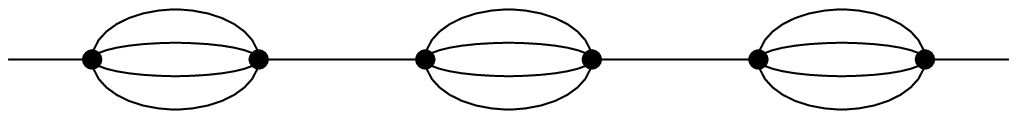}
	   & 4 & -3
	   & $\frac{1 - 9\,u^2 + 16\,u^4}{8\,u^2}$
	   & 	\includegraphics[width=1in]{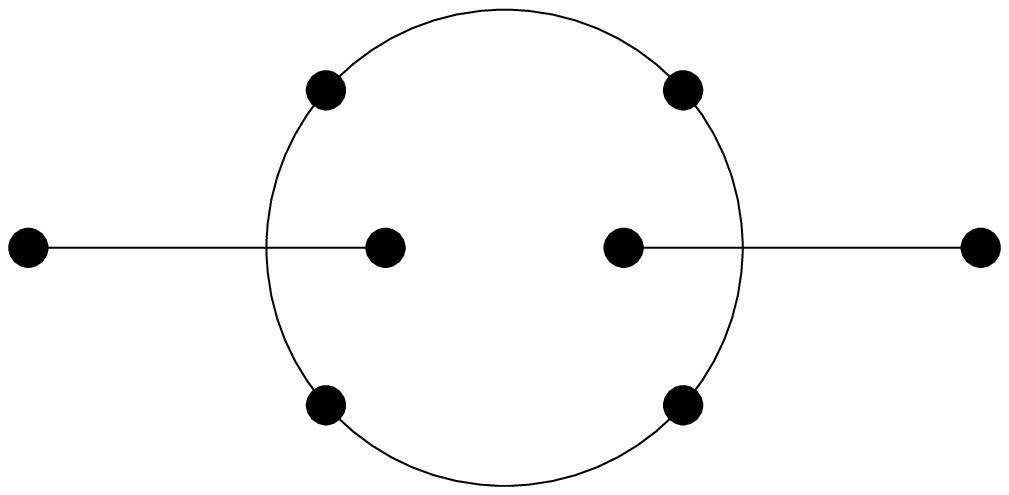}
	   & 3 \\
	   \hline
	3 & \includegraphics[width=1in]{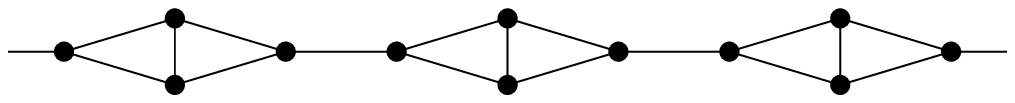}
	   & 2 & -2
	   & $\frac{1 - u + u^2 - 3\,u^3 + 2\,u^4 - 4\,u^5 + 8\,u^6}{4\,u^3}$
	   & 	\includegraphics[width=1in]{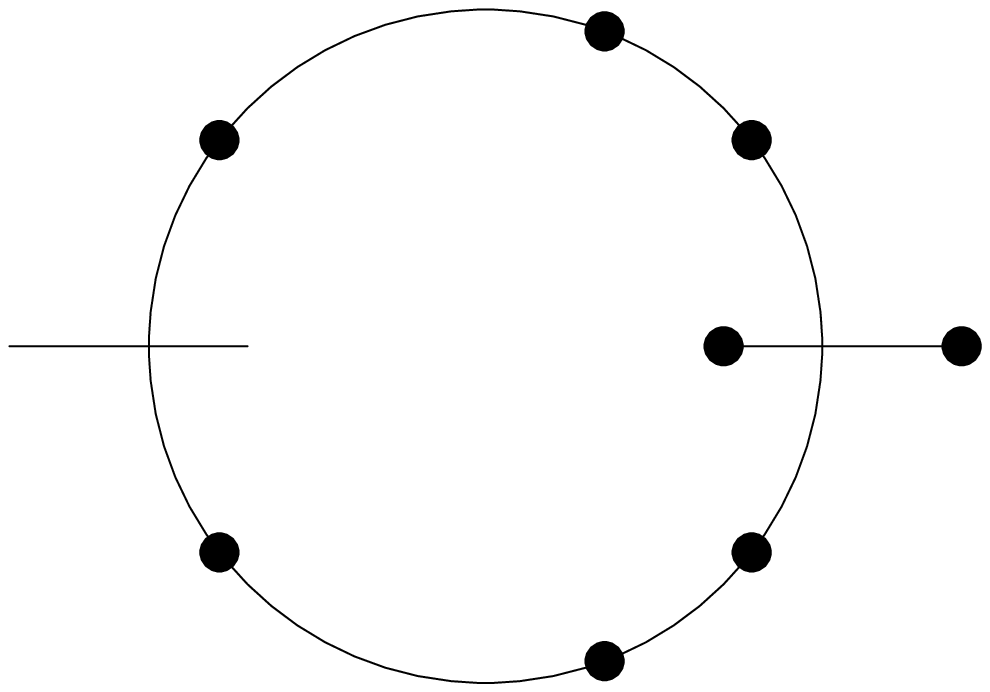}
	   & 3 \\
	   \hline
	4 &  \includegraphics[width=1in]{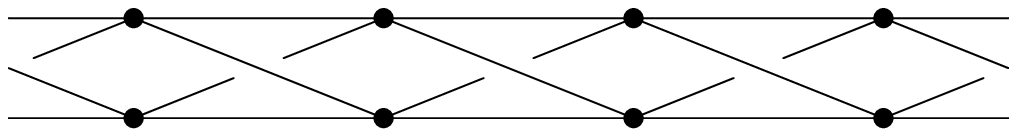}
	   & 3 & -2
	   & $\frac{1 + 3\,u^2}{4\,u}$
	   & 	\includegraphics[width=1in]{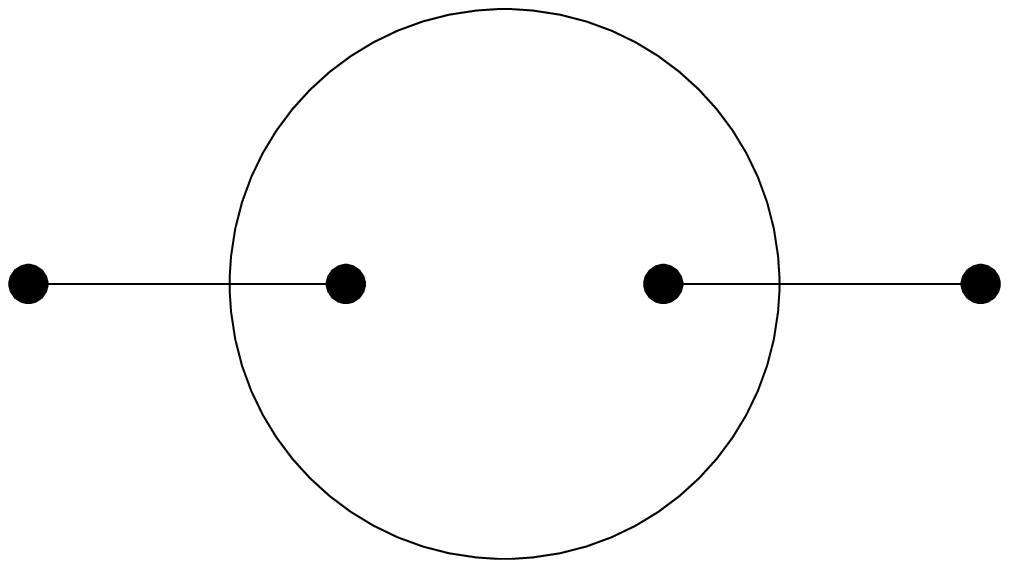}
	   & 1 \\
	   \hline
	5 & \includegraphics[width=1in]{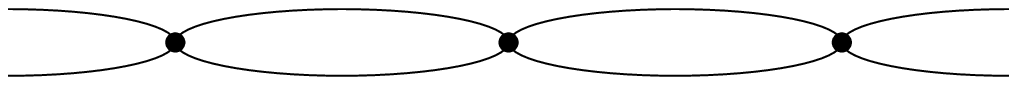}
	   & 3 & -1
 	   & $\frac{1 + 3\,u^2}{4\,u}$
	   & 	\includegraphics[width=1in]{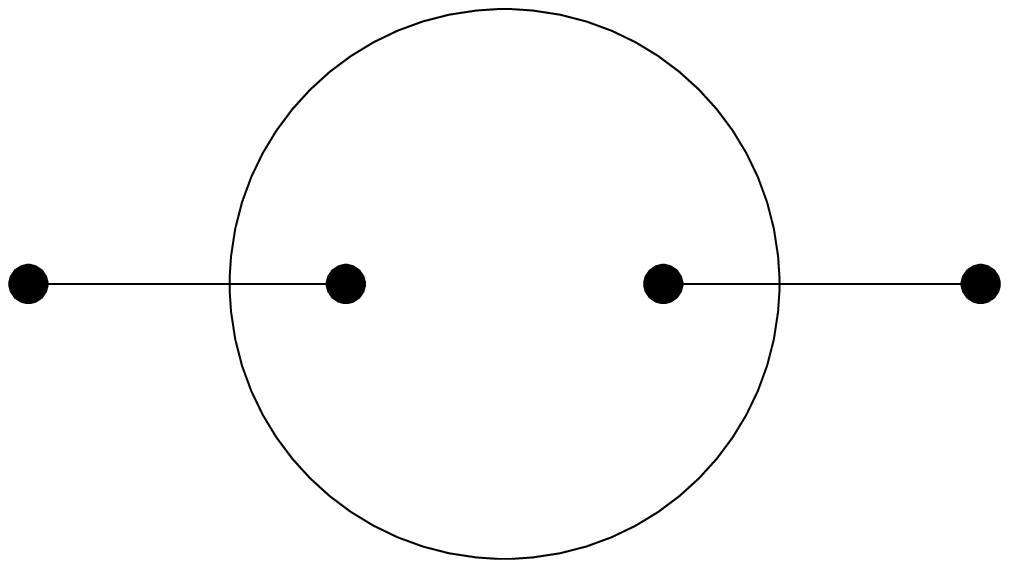}
	   & 1 \\	
	   \hline
	6 & \includegraphics[width=1in]{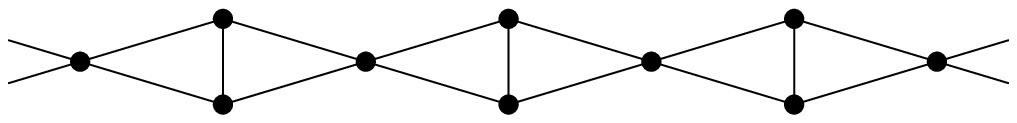}
	   & NR & -2
	   & $\frac{1 - u + u^2 - 3\,u^3 + 6\,u^4}{4\,u^2}$
	   & 	\includegraphics[width=1in]{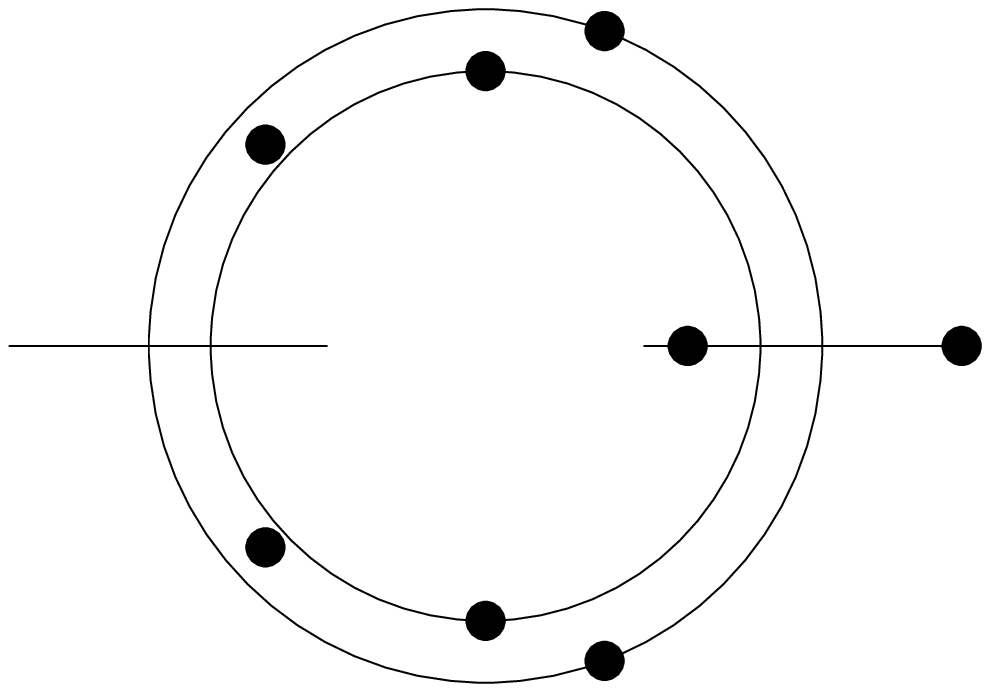}
	   & 3 \\
	   \hline
\end{tabular}
\caption{Some degree 1 Graphs}
\end{table}

\begin{Example}[A regular graph with branch point off of $C$]
	\label{ex:sawtooth}

Consider the 4-regular graph $Y$ with vertices $\Z \cup \Z$ as
shown in Figure~\ref{fig:sawtooth}.  The adjacency matrix of $Y$ is
\[
   \delta = \left(\begin{matrix}
              t + t^{-1} & 1 + t^{-1} \\
              1 + t & t+ t^{-1}
			\end{matrix}\right).
\]
$P_{u}(x)$ is degree 2, and the two roots of $P_{u}$ are
\[
   r_{\pm}(u) = \frac{1}{4u}\left(
                2+u+6u^{2}\pm \sqrt{u(4 + 9u + 12u^{2})}
				\right).
\]
From~\eqref{eq:omega}, $\Omega$ is the irreducible degree 4 
polynomial
\[
\begin{split}
   \Omega(T) = T^{4} &- \frac{1 + 4u^{2} + 9u^{4}}{u^{2}}T^{3}
             +  \frac{2 + 4\,u + 15\,u^2 + 12\,u^3 + 18\,u^4}{u^2}T^{2}\\
			 &- \frac{1 + 4u^{2} + 9u^{4}}{u^{2}}T + 1.
\end{split}
\]
The Riemann surface $S$ has four sheets covering $\C P^{1}$.
Evaluating the discriminant of $\Omega$, one has 10 points $u$ where
$Z_{Y}^{(2)}$ has duplicate values.
Checking the local behavior near those 
10 points and additionally near $u=0$, $u=\infty$, one finds that 
$S$ is unbranched at four of them.
At
\[
   u \in \left\{ 1,\frac{1}{3},\pm\frac{i}{\sqrt{3}}\right\},
\]
sheets of $S$ come together in two pairs of multiplicity two branch 
points. At
\[
   u = \left\{ 0,-\frac{9}{24} \pm 
   \frac{i}{24}\sqrt{111},\infty\right\},
\]
one pair of sheets come together in a multiplicity two branch point 
and the other two sheets are unbranched.  The pattern of branchpoints 
is shown in Figure~\ref{fig:sawtooth}, and the genus of $S$ is 3.

The most interesting thing here is that the zeta function is branched 
over 0 and $\infty$, which are not in the set $C$.
Of course, the sheet 
corresponding to the original unextended definition of $Z_{Y}^{(2)}$
is not one of the two sheets that come together at $u=0$.
The group of deck transformations of $S$ is not transitive, and the 
field extension $\C(u)[T]/(\Omega(T)) : \C(u)$ is not Galois.
\begin{figure}
	\label{fig:sawtooth}
	\includegraphics[width=1in]{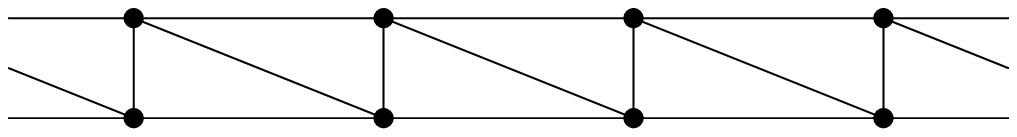}
	\includegraphics[width=1in]{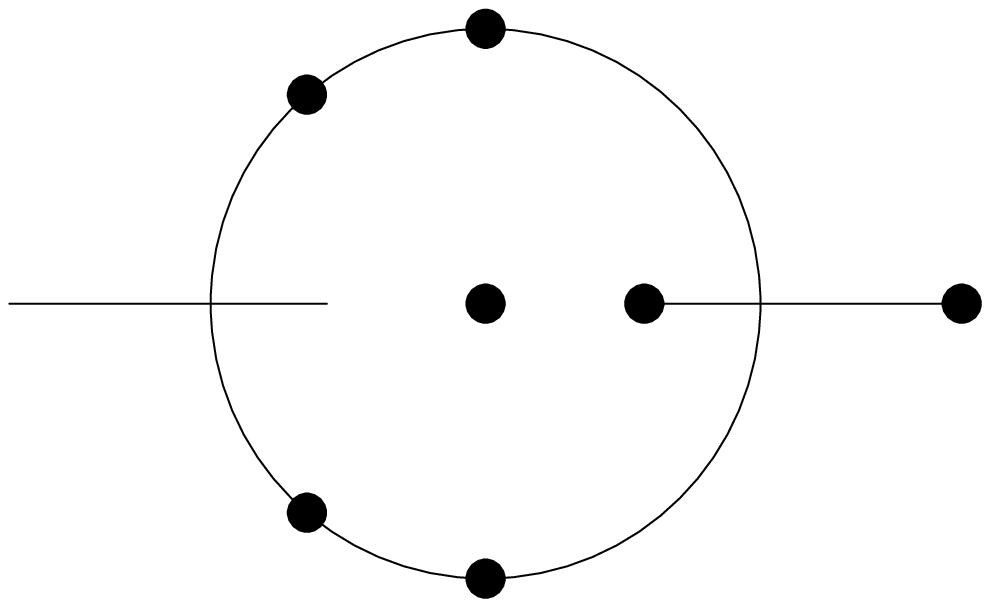}
	\caption{The graph and branchpoints for Example~\ref{ex:sawtooth}.}
\end{figure}
\end{Example}

\begin{Example}[A nontrivial reducible graph]\label{ex:triladder}
Let $Y$ have vertices $\Z \cup \Z \cup \Z$ connected as shown in
Figure~\ref{fig:triladder}.  It is the graph Cartesian
product of the line with a triangle.

Here $\Omega$ factors into a fourth degree term and the square of a 
quadratic.  The factor $\Phi$ corresponding to $Z_{Y}^{(2)}$ is the 
fourth degree term, so $S$ is four sheeted.
There are twelve branch points:
\[
	  u \in \left\{
	  \frac{1}{3},1,\frac{\pm i }{{\sqrt{3}}},
	  \frac{-3 \pm i \,{\sqrt{3}}}{6},
	  \frac{\pm 1 \pm i \,{\sqrt{11}}}{6},
	  -\frac{1}{4}  + \frac{{\sqrt{\frac{7}{3}}}}{4}
	      \pm \frac{i}{2}\,{\sqrt{\frac{1}{2} + 
	      \frac{{\sqrt{\frac{7}{3}}}}{2}}}
	  \right\},
\]
shown in Figure~\ref{fig:triladder}.
At each $u$, sheets come together in two pairs of multiplicity two 
branch points, so the genus of $S$ is 9.

\begin{figure}\label{fig:triladder}
	\includegraphics[width=1in]{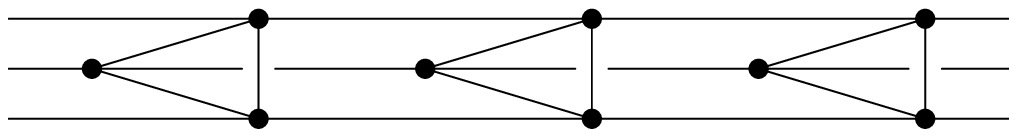}
	\includegraphics[width=1in]{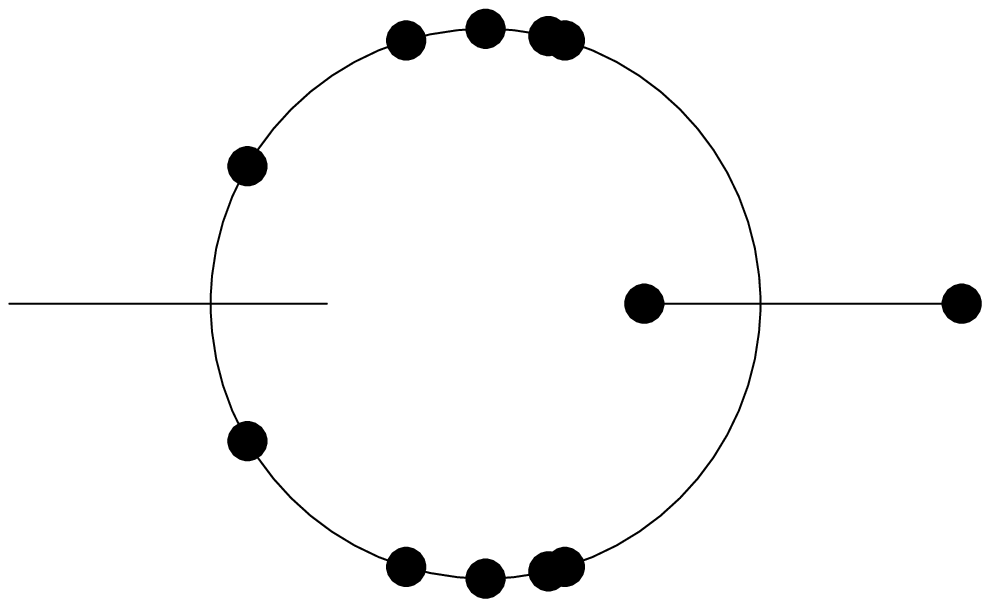}
	\caption{The graph and branchpoints for Example~\ref{ex:triladder}.}
\end{figure}

Even though $\Omega$ is reducible, all branch points still lie on the set 
$C$.  Here $Z_{Y}^{(2)}$ must still satisfy the functional 
equation~\eqref{eq:func}, because the involution $\tilde{\tau}$ 
preserves $S$ - the other two irreducible factors of $\Omega$
are degree 2 while $\Phi$ is degree 4.
\end{Example}

%
%
\bibliographystyle{plain}
\bibliography{zetaz}

\end{document}